\documentclass[12pt]{article}

\usepackage{amscd}

\begin{document}

\title{Compatible systems of mod $p$ Galois representations: II 
}

\author{Chandrashekhar Khare}

\date{}

\maketitle

\newtheorem{theorem}{Theorem}
\newtheorem{lemma}{Lemma}
\newtheorem{prop}{Proposition}
\newtheorem{cor}{Corollary}
\newtheorem{example}{Example}
\newtheorem{conjecture}{Conjecture}
\newtheorem{definition}{Definition}
\newtheorem{quest}{Question}
\newtheorem{memo}{Memo}
\newcommand{\rhobar}{\overline{\rho}}
\newcommand{\Sha}{{\rm III}}
\newtheorem{conj}{Conjecture}

\section{Introduction}

Compatible systems of $n$-dimensional, 
mod $p$ representations of absolute Galois groups 
of number fields were considered by Serre
in his study of openness of 
images of adelic Galois representations arising from elliptic curves 
in [S2]. In [K] the author considered 
{\it abstract} compatible systems of $n$-dimensional, 
mod $p$ representations of 
the absolute Galois group of ${\bf Q}$ and determined them in
the one-dimensional case.

\subsection{Definition and main theorem}

We begin with the definition of compatible systems of $n$-dimensional, 
mod $p$ representations 
of the absolute Galois group of a number field $K$.

\begin{definition}\label{compatible}

\begin{itemize} \item Let $K$ and $L$ be number fields and $S,T$ finite sets of places 
 of $K$ and $L$ respectively. An $L$-rational (resp., $L$-integral) 
 strictly compatible system $\{\rho_{\wp}\}$ 
 of $n$-dimensional mod $\wp$ representations of
 $G_K:={\rm Gal}(\overline{K}/K)$ with defect set $T$ and ramification set $S$,
 consists of giving for each finite place $\wp$ of $L$ not in $T$ a 
 continuous, semisimple representation $$\rho_{\wp}:G_K \rightarrow GL_n({\bf
 F}_{\wp}),$$ for ${\bf F}_{\wp}$ the residue field of ${\cal O}_L$ at
 $\wp$ of characteristic $p$, that is \begin{itemize} \item  unramified at the places
 outside $S \cup$ $\{$ places of $K$
 above $p$ $\}$ \item there is a monic polynomial $f_r(X) \in L[X]$ (resp., 
 $f_r(X) \in {\cal O}_L[X]$) such that  for each place 
 $r$ of $K$ not in $S$ and for all places $\wp$ of
 $L$ not in $T$, coprime to the residue characteristic of $r$, and
 such that $f_r(X)$ has coefficients that are integral at $\wp$, 
 the characteristic polynomial of $\rho_{\wp}({\rm Frob}_r)$
 is the reduction of $f_r(X)$ mod $\wp$, where ${\rm Frob}_r$ is the
 conjugacy class of the Frobenius at $r$ in the Galois
 group of the extension of $K$ that is the fixed field of the kernel of
  $\rho_{\wp}$. \end{itemize} 
\item If the prime to $p$ part of the Artin conductor of
 $\rho_{\wp}$ 
 is bounded independently of $\wp$ 
 we say that the system $\{\rho_{\wp}\}$ has
 bounded conductor.\end{itemize}
\end{definition}
In the course of the paper we will often suppress the sets $S,T$ from the notation. It might be prudent
to impose the condition of 
integrality of the roots of $f_r(X)$ outside primes coprime to the
residue characteristic of primes in $T$ and $r$ in Definition
\ref{compatible}, as only then are the $\rho_{\wp}$'s determined
a priori by the compatible system data $\{f_r(X)\}$. But as we do not need this
in the main result of the paper, 
Theorem \ref{main} below, we stick to our less stringent requirements.

One would like to prove that
a strictly compatible system 
arises {\it motivically}, i.e., ``from the mod $p$
\'etale cohomology of a variety $X_{/K}$ as $p$ varies''.
In the case of one-dimensional
strictly compatible systems 
this is interpreted as saying that it
arises from a Hecke character.
It is only in the one-dimensional case thanks to class field theory
that at the moment one has a realistic chance of describing compatible
sytems in any degree of generality. We do this below and prove 
in Section 4.2 the
following theorem that is the main result of the paper (for
any unexplained terms see Section 4.1, in particular 
Definition \ref{compatible2}). 

\begin{theorem}\label{main}
 A $L$-rational strictly compatible system $\{\rho_{\wp}\}$ 
 of one-dimensional mod $\wp$ representations of
 ${\rm Gal}(\overline{K}/K)$ arises from a Hecke character.
\end{theorem}

\begin{cor}\label{taniyama}
  An $L$-rational strictly compatible system of one-dimensional 
  $\wp$-adic representations
  $\{\rho_{\wp,\infty}\}$ as in I-11 of [S]  
  arises from a Hecke character.
\end{cor}

\noindent{\bf Proof of corollary:} This follows from the general observation
that, given a strictly compatible system of $\wp$-adic 
representations as in I-11 of [S], we can reduce it mod $\wp$ 
and semisimplify to get a rational strictly compatible system in the
sense  of Definition \ref{compatible} such that $f_r(X)$ is integral away from the primes
that have the same residue characteristic as $r$. Thus the resulting strictly compatible 
system of mod $\wp$ representations determines the 
system of $\wp$-adic representations. From 
this observation and Theorem \ref{main} 
the corollary follows.

\vspace{3mm}

The corollary is Proposition 1.4 in [Sch]. It is
deduced there as a consequence of a deep result of Waldschmidt in
transcendental number theory that as a corollary proves the much
stronger result that 
even a {\it single}, algebraic, one-dimensional $\wp$-adic representation
arises from a  Hecke character (see [H]).
The corollary in a weaker form, assuming a supplementary 
purity hypothesis, is also  
an old theorem of Taniyama (see Theorem 1 of  [Tani]) that he proved using methods 
quite different from the present paper.  But both in [Tani] and in
the present work the compatibility hypotheses 
are used to describe 1-dimensional
compatible systems using purely algebraic methods, and 
Theorem \ref{main} can be addressed using only such methods.

We point out a corollary that is immediate from the theorem
but that sets us up to formulate two conjectures
in the higher dimensional case.

\begin{cor}
 An $L$-rational strictly compatible system $\{\rho_{\wp}\}$ 
 of one-dimensional mod $\wp$ representations of
 ${\rm Gal}(\overline{K}/K)$ 
 lifts to a strictly compatible 
 system of $\wp$-adic representations. Further it
 is of bounded conductor, and satisfies the
 purity and integrality properties of Conjecture \ref{galois} below.
\end{cor}

\subsection{Conjectures and their discussion}

Guided by the results above we make the following two conjectures. In the first we propose a reciprocity law for compatible
systems of mod $\wp$ representations as the following ``meta-conjecture'':

\begin{conjecture}\label{meta}
  Any strictly compatible $L$-integral 
  system as in Definition \ref{compatible} arises
  motivically.
\end{conjecture}

To be more specific we propose the following 
purely Galois theoretic conjecture:

\begin{conjecture}\label{galois} Let $\{\rho_{\wp}\}$ be a strictly compatible
system as in Definition \ref{compatible}.
\begin{itemize}\item (Lifting)
  It lifts to  (i.e., is the reduction up to semisimplification of)
  a strictly compatible system of semisimple $\wp$-adic representations.
\item (Bounded conductor) It is of
  bounded conductor: more precisely the prime-to-$p$ part of the Artin
  conductor of $\rho_{\wp}$ is independent of $\wp$ for almost all $\wp$. 
\item (Purity) Assume that $\rho_{\wp}$ is irreducible for almost all $\wp$. Then the roots of $f_r(X)$ for primes $r$ not in $S$ are 
of absolute value $|{\rm
Nm}(r)|^t$ with respect to all embeddings of $\overline{{\bf
Q}}$ in ${\bf C}$, and for an integer or 
half-integer $t$ that is idependent 
of $r$. Here ${\rm Nm}$ is the norm map to ${\bf Q}$.
\item (Integrality) $\{\rho_{\wp} \otimes {\rm Nm}_{\wp}^m\}$ is integral where ${\rm
  Nm}_{\wp}$ is the $\wp$-adic cyclotomic character for some integer
$m$. \end{itemize}
\end{conjecture}

Because of known properties of Galois representations
which arise from geometry one expects that Conjecture \ref{meta} implies 
Conjecture \ref{galois} in the integral case. 
Note that the first part of Conjecture \ref{galois} implies
integrality of the roots of $f_r(X)$ at primes coprime to the
residue characteristic of primes in $T$ and $r$, while the 
last part of the conjecture implies
integrality of the roots of $f_r(X)$ at primes coprime to 
the residue characteristic of $r$. 
Haruzo Hida has also pointed out that one might expect that the minimal
field of rationality of a strictly compatible system as in Definition
\ref{compatible} is totally real or CM: this is again verified in
the one-dimensional case by Theorem \ref{main}. Our conjectures may be
regarded
as an analog for compatible mod $p$ systems of Galois representations
of the well-known conjectures of Fontaine and Mazur, cf. [FM].
Proving that a
strictly compatible system $\{\rho_{\wp}\}$ arises motivically 
is much stronger than proving that an individual $\rho_{\wp}$
arises motivically: this is in contrast to what happens in the case of compatible
systems of $\wp$-adic representations. Nevertheless in the mod $\wp$ case too, there is a link between reciprocity conjectures 
about the individual representations $\rho_{\wp}$ being motivic and
the compatible system arising motivically, if we grant the properties
of the compatible system $\{\rho_{\wp}\}$ of Conjecture \ref{galois}.
The reader is invited to consult [S1] where such a link is established
which allows Serre to deduce the Shimura-Taniyama conjecture as
a consequence of the conjectures in [S1].

Theorem \ref{main} proves these conjectures
for 1-dimensional strictly compatible systems. On the other hand there are
pairs of 1-dimensional mod $p$ and mod $q$ representations
that do not arise simultaneously from a Hecke character. In [KK] it is determined when 
a given pair of one-dimensional mod $p$ and mod $q$ representations 
arises from a Hecke character.

Here is the plan of the paper: In Section 2 we prove a descent result (Proposition
\ref{lift} below) related to the lifting property.
In Section 3 we generalise a result of [CS] that is needed for
the proof of Theorem \ref{main}. The proof is carried out in Section 4.
In Section 5 we partially generalise the results of Section 4 to
abelian semisimple compatible systems. Although we do not have a complete result 
in this case, as there is an essential difficulty that we are unable to overcome,
the methods used in this section may be of independent interest.
In particular we would like to draw attention to an analog of Artin's conjecture on primitive roots,
for the ring that arises from restriction of scalars to ${\bf Z}$ of the ring of integers of a number field,
that we formulate in this section. 

\section{Descent results}

We prove a descent result in the context of the lifting statement of 
Conjecture \ref{galois}. We first define what we mean by {\it weakly
compatible systems}:

\begin{definition}\label{weak} \begin{itemize} \item
 Let $K$ and $L$ be number fields and $S,T$ finite sets of places 
 of $K$ and $L$ respectively. A $L$-rational (resp., $L$-integral) 
 weakly compatible system $\{\rho_{\wp,\infty}\}$ 
 of $n$-dimensional $\wp$-adic representations of
 ${\rm Gal}(\overline{K}/K)$ with defect $T$ and ramification set $S$,
 consists of giving for each finite place $\wp$ of $L$ not in $T$ a 
 continuous semsimple representation $$\rho_{\wp,\infty}:G_K \rightarrow
 GL_n(L_{\wp}),$$ for $L_{\wp}$ the completion of $L$ at
 $\wp$ whose residue field has characteristic $p$, 
 that is \begin{itemize} \item  unramified at the primes 
 outside $S \cup$ $\{$ places of $K$
 above $p$ $\}$ \item for a place $r$ of $K$ not in $S$ 
 there is a monic polynomial $f_r(X) \in L[X]$ (resp., 
 $f_r(X) \in {\cal O}_L[X]$) such that for almost all places $\wp$ of
 $L$, the characteristic polynomial of $\rho_{\wp,\infty}({\rm Frob}_r)$
 is $f_r(X)$, where ${\rm Frob}_r$ is the
 conjugacy class of the Frobenius at $r$ in the Galois
 group of the extension of $K$ that is the fixed field of the kernel of
  $\rho_{\wp,\infty}$. \end{itemize} 
 \item A $L$-rational (resp., $L$-integral) 
 weakly compatible system $\{\rho_{\wp}\}$ 
 of $n$-dimensional mod $\wp$ representations of
 ${\rm Gal}(\overline{K}/K)$ with defect $T$ and ramification set $S$,
 consists of giving for each finite place $\wp$ of $L$ not in $T$ a 
 continuous semsimple representation $$\rho_{\wp}:G_K \rightarrow GL_n({\bf
 F}_{\wp}),$$ for ${\bf F}_{\wp}$ the residue field of ${\cal O}_L$ at
 $\wp$ of characteristic $p$, that is \begin{itemize} \item  unramified at the primes 
 outside $S \cup$ $\{$ places of $K$
 above $p$ $\}$ \item for a prime $r$ of $K$ not in $S$ 
 there is a monic polynomial $f_r(X) \in L[X]$ (resp., 
 $f_r(X) \in {\cal O}_L[X]$) such that for almost all placess $\wp$ of
 $L$, the characteristic polynomial of $\rho_{\wp}({\rm Frob}_r)$
 is the reduction of $f_r(X)$ mod $\wp$, where ${\rm Frob}_r$ is the
 conjugacy class of the Frobenius at $r$ in the Galois
 group of the extension of $K$ that is the fixed field of the kernel of
  $\rho_{\wp}$. \end{itemize} \end{itemize}
\end{definition}

We do not know if
there are weakly compatible systems in the mod
$\wp$ setting that are not strictly compatible for some large but
finite defect set. (Note that we might use the equivalent
words place and prime below.)

\begin{prop}\label{lift}
 Let $K'$ be a finite Galois extension of $K$. Consider a  
 $L$-integral strictly compatible system $\{\rho_{\wp}\}$ 
 of mod $\wp$ representations of
 $G_K:={\rm Gal}(\overline{K}/K)$, with the further property that
 $\rho_{\wp}|_{G_{K'}}$ is absolutely irreducible for almost all
 $\wp$. Then if the strictly compatible system
 $\{\rho_{\wp}|_{G_{K'}}\}$ lifts to a strictly compatible $L'$-integral 
 $\wp$-adic system $\{\rho_{\wp,K',\infty}\}$ for some number field $L'$
 with finite defect and exceptional sets,
 then the system $\{\rho_{\wp}\}$ lifts to a weakly compatible $L$-integral 
 $\wp$-adic system $\{\rho_{\wp,K,\infty}\}$.
\end{prop}

\noindent{\bf Proof:}  We only consider the $\rho_{\wp}$'s for $\wp$'s
such that their residue characteristic is prime to the degree
$[K':K]$ and for which $\rho_{\wp}|_{G_{K'}}$ is irreducible. 
We  assume without loss of generality
that $L'$ contains $L$ and that we are considering $L'$-rational
compatible systems and denote the places of $L'$ by $\wp$ etc.
Then: \begin{itemize} \item 
$\rho_{\wp}|_{G_{K'}}$ lifts to $\rho_{\wp,K',n}$ 
the mod $\wp^n$ representation that is the reduction of $\rho_{\wp,K',\infty}$,
\item Note that the residue characteristic of $\wp$ 
is prime to the degree
$[K':K]$, and $\rho_{\wp}|_{G_{K'}}$ extends to the representation
$\rho_{\wp}$ of $G_K$. It is easy to see from this, 
by computing cohomological obstructions (see Proposition 1.1 of [Clo]
for a similar argument),
that for each $n$, the representation $\rho_{\wp,K',n}$ extends to a representation 
$\rho_{\wp,K,n}$ of $G_K$, using that
$\rho_{\wp,K',n}$ satisfies descent
data as $\{\rho_{\wp,K',\infty}\}$ is a strictly compatible system
and $\{\rho_{\wp,K'}\}$ comes by restriction from $\{\rho_{\wp}\}$.
Then we are done by invoking Carayol's theorem, see [Ca], as in [Clo].  
\item $\rho_{\wp}|_{G_{K'}}$ being 
irreducible, any extension of $\rho_{\wp}|_{G_K'}$ to $G_K$ is unique
up to twisting by characters of ${\rm Gal}(K'/K)$, and thus we may
assume that for each $n$, $\rho_{\wp,K,n}$ reduces mod $\wp$ to
$\rho_{\wp}$. \end{itemize} Thus for all but finitely many $\wp$, 
$\rho_{\wp}$ lifts to a $\wp$-adic
representation $\rho_{\wp,K,\infty}$ whose restriction to $G_{K'}$ 
is $\rho_{\wp,K',\infty}$. Now we claim that the
$\{\rho_{\wp,K,\infty}\}$'s for such $\wp$'s 
form the desired weakly compatible lift of the
$\rho_{\wp}$'s. First observe that  $\rho_{\wp,K,\infty}$ is
unramified outside primes above $p$ and a fixed finite set $S'$ 
that does not depend on $p$, that we see using the fact that $K'/K$ is ramified at only finitely many primes. 
Then observe that for almost all primes $r$
of $K$ there are only finitely many possibilities 
for the roots of the characteristic polynomial of ${\rm Frob}_r$
in $\rho_{\wp,K,\infty}$ as $\wp$ varies, and $r$ is fixed.
This is because $\{\rho_{\wp,K',\infty}\}$
is in particular a weakly compatible system. Thus using the
srictly compatible system $\{\rho_{\wp}\}$ it follows that
the characteristic polynomial of $\rho_{\wp,K,\infty}({\rm Frob}_r)$
is $f_r(X)$ for almost all primes $\wp$ and fixed $r$. 

\vspace{3mm}

\noindent{\bf Remark:} We unfortunately do not know how to prove 
the more desirable result that
under the conditions above  
$\{\rho_{\wp}\}$ lifts to a strictly compatible system of $\wp$-adic
representations, although the weakly compatible $\wp$-adic system
constructed above {\it should be} strictly compatible
for a suitably large defect set. This is because
the compatible system $\{\rho_{\wp}\}$ serves as rigidifying data in
the sense that for almost all primes, more precisely all primes
of residue characteristic prime to $[K':K]$, the extension 
$\rho_{\wp,K,\infty}$ of $\rho_{\wp,K',\infty}$ is 
{\bf determined} by the requirement that it reduces
to $\rho_{\wp}$: this follows from the irreducibility hypothesis
on $\{\rho_{\wp}|_{G_{K'}}\}$. 

\vspace{3mm}

We now prove a more specific descent result in the context of 
Conjecture \ref{meta} for 1-dimensional representations that will
be useful in the proof of Theorem \ref{main} (see Definition
\ref{compatible2} of Section 4.1 below for unexplained terms).

\begin{lemma}\label{util}
 Let $K'$ be a finite extension of $K$. Then a 
 $L$-rational strictly compatible system $\{\rho_{\wp}\}$ 
 of 1-dimensional mod $\wp$ representations of
 $G_K:={\rm Gal}(\overline{K}/K)$ with defect $T$ and exceptional set $S$
 and with bounded conductor arises from a Hecke character if and only
 if the strictly compatible system $\{\rho_{\wp}|_{G_{K'}}\}$ arises from a
 Hecke character.
\end{lemma}

\noindent{\bf Proof:} Only one direction needs a proof. We may assume
without loss of generality that $K'$ is a Galois extension of ${\bf
Q}$. So assume that the strictly compatible system 
$\{\rho_{\wp}|_{G_{K'}}\}$ arises from a Hecke character $\chi$
of the idele group of $K'$ and of infinity
type $(m_{\sigma})_{\sigma \in {\rm Gal}(K'/{\bf Q})}$. Using that the strictly compatible system
$\{\rho_{\wp}|_{G_{K'}}\}$ arises by restriction from a strictly compatible
system of $G_K$ we see that $\chi^{\sigma}=\chi$ for all $\sigma \in
{\rm Gal}(K'/K)$. Consider a principal prime ideal $(a)$ of $K'$
that lies above a prime of ${\bf Q}$ that splits completely in
$K'$. Then as $\chi(\sigma(a))=\chi(a)$ for all $\sigma \in
{\rm Gal}(K'/K)$ we deduce that $m_{\sigma}=m_{\sigma'}$ whenever
the restrictions of $\sigma,\sigma' \in {\rm Gal}(K'/{\bf Q})$ to $K$ are equal, and thus for
each embedding $\sigma$ of $K$ we can put without ambiguity $m_{\sigma}$
equal to $m_{\sigma'}$ for any $\sigma' \in {\rm Gal}(K/{\bf Q})$ that
extends $\sigma$. Now consider the 
algebraic character $K^* \rightarrow {\bf C}^*$ given by
$k \rightarrow \Pi_{\sigma} \sigma(k)^{m_{\sigma}}$ where $\sigma$ runs through
embeddings of $K$ in ${\bf C}$. It is easy to see that this character is trivial on a subgroup of
finite index of the units ${\cal O}_K^*$. Thus it is trivial
by Th\'eor\`eme 1 of [C] on units congruent to 1 mod $\sf n$ of ${\cal
O}_K^*$ for some ideal $\sf n$ of ${\cal O}_K$.
It follows that there is a Hecke character $\chi'$ 
for $K$ such that $\rho_{\wp} \otimes \tilde{\chi'}_{\wp}^{-1}$ where
$\{ {\chi'}_{\wp} \}$ is the compatible mod $\wp$-system that $\chi'$
gives rise to as described in the next section factors through a
fixed finite extension of $K$. Observe that strictly
compatible systems as in Definition \ref{compatible} that factor through
the Galois group ${\rm Gal}(K''/K)$ of a fixed finite extension $K''$
of $K$ arise as reductions mod primes 
of a continuous 
representation of $G_K$ into $GL_n({\cal O})$ with finite image: here ${\cal O}$ the ring of integers of a number field. This finishes
the proof.

\vspace{3mm}

We say that a compatible system $\{\rho_{\wp}\}$
as in Definition \ref{compatible} arises from an Artin representation
if it arises
from reducing a 
continuous 
representation $\rho:G_K \rightarrow GL_n({\cal O})$ with finite image,
where ${\cal O}$ the ring of integers of a number field $L$,
modulo the primes of $L$ and semisimplifying.
We quote a result in Section 8 of 
Deligne-Serre (cf. [DS]: we thank C.S. Rajan for this reference)
that characterises compatible systems $\{\rho_{\wp}\}$
that arise from Artin representations and refines the observation
towards the end of the proof above.

\begin{prop}\label{DS}(Deligne-Serre)
  Let $\{\rho_{\wp}\}$ be a $L$-integral compatible system
  where we allow the defect set $T$
  to be any set whose complement 
  in the set of places of $L$ is infinite. 
  If $|{\rm im}(\rho_{\wp})|$
  is bounded independently of $\wp$, then $\{\rho_{\wp}\}$
  arises from an Artin representation.
\end{prop}

\vspace{3mm}

\noindent{\bf Remarks:} \begin{enumerate} \item One can ask for another subtler 
characterisation of integral strictly compatible systems
$\{\rho_{\wp}\}$ that arise from an Artin representation as those that are
unramified outside a fixed
finite set of places that is independent of $\wp$.
\item Compatible systems arising form Artin representations were
used to prove that Serre's conjectures in [S1] imply the
modularity of 2-dimensional irreducible odd complex representations
of $G_{\bf Q}$ in [K1].
\end{enumerate}

\section{A result of Corrales and Schoof}

We need the following straightforward generalisation of Theorem 1 of [CS].

\begin{prop}\label{schoof}
Let $a_1,\cdots,a_n,c$ be a finite set of elements of $K^*$ and $\ell$ a
rational prime. If for almost all primes
$\wp$ of $K$ the subgroup generated by the image of $a_i$'s in
$({\cal O}_k/\wp)^*$ contains $c^{t_{\wp}}$ for some integer $t_{\wp}$ that
depends on $\wp$ and is prime to $\ell$,
then $c^t=\Pi_{i=1}^{i=n}a_i^{m_i}$ for some integer $t$ with $t$ prime to $\ell$, and 
for some integers $m_i$, $i=1,\cdots,n$.
If further $t_{\wp}$ can be taken to be 1 for almost all $\wp$, then 
$c=\Pi_{i=1}^{i=n}a_i^{m_i}$ for some integers $m_i$, $i=1,\cdots,n$.
\end{prop}

\noindent{\bf Proof:} Without loss of generality we can and will assume $\sqrt{-1} \in K$.
The proof is basically the same as that of Theorem 1 of
[CS]: we will briefly sketch the proof using the same breakdown into
steps as in [CS].\begin{itemize} \item Step 1: Let $q$ be a power of
$\ell$. Consider $K(\zeta_q,a_1^{1/q},\cdots,a_n^{1/q})$ and
$K(\zeta_q,c^{1/q})$. Observe that for all but finitely many
exceptions
a prime $\wp$ of $K$ that lies above a prime of ${\bf Q}$ that splits
completely
in $K/{\bf Q}$ splits in $K(\zeta_q,c^{1/q})$ (resp.,
$K(\zeta_q,a_1^{1/q},\cdots,a_n^{1/q})$) if and only if $p$ (the
characteristic of the residue field at $\wp$) is 1 mod
$q$, and $c$, or equivalently $c^{t_{\wp}}$ as $t_{\wp}$ is prime to
$\ell$, is a $q$th power in ${\bf F}_{\wp}^*$ (resp., $a_1,\cdots,a_n$
are $q$th powers in ${\bf F}_{\wp}^*$). Thus by assumptions
of the proposition almost all split primes of $K$ that split
in $K(\zeta_q,a_1^{1/q},\cdots,a_n^{1/q})$
also split in
$K(\zeta_q,c^{t_{\wp}/q})=K(\zeta_q,c^{1/q})$. Thus by 
the Frobenius density theorem we conclude that $K(\zeta_q,c^{1/q}) \subset
K(\zeta_q,a_1^{1/q},\cdots,a_n^{1/q})$.
\item Step 2: Using Kummer theory and duality for finite abelian
groups we deduce from Step 1 that the subgroup generated by the images of $a_i$'s 
in $K(\zeta_q)^*/K(\zeta_q)^{*q}$ contains the image of $c$. From 
Lemma 2.1 (ii) of [CR] (see also remark on page 39 of [C]), 
we have that the natural map 
$K^*/K^{*q} \rightarrow K(\zeta_q)^*/K(\zeta_q)^{*q}$ is injective,
and thus we deduce that the
subgroup generated by the images of $a_i$'s 
in $K^*/K^{*q}$ contains the image of $c$. \item Step 3: Consider 
$A:={\cal O}_{K,T'}^*/\langle a_1,\cdots,a_n \rangle$ where
${\cal O}_{K,T'}^*$ are the $T'$-units with $T'$ consisting of all
places dividing any of the $a_i$'s or $c$, the infinite places and
the places where the hypotheses of  the statement of the 
proposition are not satisfied. 
We know from Step 2 that $c$ is in $A^q$ for all $q$'s
that are powers of $\ell$. Using the fact that ${\cal O}_{K,T'}^*$ is finitely
generated (Dirichlet unit theorem), we conclude
that the image of $c$ in $A$ is torsion of order prime to $\ell$ and thus we are done.
\end{itemize}

The last line of the proposition follows from the fact that assuming
the $t_{\wp}$'s to be 1 we can work with all prime powers $q$.

\section{Reciprocity for one-dimensional compatible systems}

\subsection{Generalities about Hecke characters}

We recall the definition of Hecke characters and the association of
compatible systems of $p$-adic representations to them (see Chapter II of
[S], or [W] and [Tani]). There 
are equivalent ways of looking at Hecke characters that we will
need below and we briefly recall these.
We index as usual the real (resp., complex) places of $K$ by embeddings
$\sigma$ (resp.,
pairs of elements $\{\sigma,c\sigma\}$ where $c$ is complex
conjugation) of $K$ in ${\bf C}$. 
Let $I$ be the group of ideles of $K$ and ${(K_{\infty}^{\times})}^0$
be the connected component of the identity of the product of the
completions of $K$ at the archimedean places.

\begin{definition}\label{hecke}
A Hecke character $\chi$ is a 
continuous homomorphism 
$\chi: I/K^*  \rightarrow {\bf C}^*$
such that $$\chi|_{{(K_{\infty}^{\times})}^0}(x)=
\Pi_{\sigma \ real}x_{\sigma}^{n_{\sigma}}  \Pi_{\sigma \ complex}
x_{\sigma}^{n_{\sigma}} \overline{x}_{\sigma}^{n_{c\sigma}}$$ for
integers $n_{\sigma},n_{c\sigma}$ and with
$x_{\sigma}$ the components of $x$. We say that the tuple of
integers $(n_{\sigma})_{\sigma}$ is the infinity type of $\chi$, and
say that $\chi$ is unramified at a finite place $v$ if the units $U_v$
at $v$ are in the kernel of $\chi$. The conductor of $\chi$ is the
largest ideal $\sf n$ such that elements of the finite ideles $I^{(\infty)}$ congruent to
1 mod $\sf n$ are in the kernel of $\chi$.
\end{definition}

From $\chi$ we get a continuous homomorphism
$\chi_0: I/ {(K_{\infty}^{\times})}^0 \rightarrow {\bf C}^*$
defined by $\chi_0(x)=\chi(x)\Pi_{\sigma \
real}x_{\sigma}^{-n_{\sigma}}  
\Pi_{\sigma \ complex}
x_{\sigma}^{-n_{\sigma}} \overline{x}_{\sigma}^{-n_{c\sigma}}$ whose kernel
is open and takes values in a
sufficiently large subfield $L$ of ${\bf C}$ 
which is a finite extension of ${\bf Q}$. 

For any finite place
$\wp$ of $L$ above a prime $p$ of ${\bf Q}$ the morphism $\psi:K^*
\rightarrow L_{\wp}^*$, $\psi(x)=\Pi_{\sigma}\sigma(x)^{n_{\sigma}}$,
extends to a continuous morphism $\psi_{p}:(K \otimes {\bf Q}_p)^*
\rightarrow L_{\wp}^*$, and we define $\chi_{\wp}:I
/K^*{(K_{\infty}^{\times})}^0 \rightarrow L_{\wp}^*$ by $\chi_{\wp}=\chi_0
.(\psi_{p}\alpha_p)$ where $\alpha_p$ is the projection of $I$ to the
components at places above $p$. The kernel of $\chi_{\wp}$ is open in
$I^{(p)}$ the ideles concentrated at places away from those dividing
$p$. Using the isomorphism of class field
$G_K^{ab} \simeq I
/K^*{(K_{\infty}^{\times})}^0$, we see 
that $\{\chi_{\wp}\}$ forms a compatible
system of 1-dimensional $\wp$-adic representations of $G_K$ in a
natural way. Since $G_K$ is compact, $\chi_{\wp}$
takes values in the units ${\cal O}_{\wp}^*$ of $L_{\wp}^*$ and
thus can be reduced mod $\wp$.

\begin{definition}\label{compatible2}
We say that a strictly compatible system of (one-dimensional) mod $\wp$
representations $(\rho_{\wp})$ as in Definition \ref{compatible}
arises from a Hecke character $\chi$ if $\rho_{\wp}=\tilde{\chi_{\wp}}$
where $\tilde{\chi_{\wp}}$ is the reduction of $\chi_{\wp}$ mod $\wp$
at all primes $\wp$ not in $T$.
\end{definition}

For a fractional ideal $\sf n$ of ${\cal O}_K$, which we identify with a 
sequence of integers $(m_v)$ for $v$ 
running through finite places of $K$ and $m_v=0$ for almost all $v$, 
define the subgroup
$U_{\sf n}$ of the ideles $I$ of $K$ to be the product $U_{v,\sf n}$
where $U_{v,\sf n}$ is the connected component of $K_v^*$ if
$v$ is an infinite place, and the units $U_v$ congruent to 1 modulo
the $m_v$ th power of the maximal ideal if $v$ is a finite place. Thus
$K^* \cap U_{\sf n}$ are the totally positive units $E_{\sf n}$ 
of ${\cal O}_K$ congruent to 1 mod $\sf
n$. Let $I_{\sf n}$ be the quotient $I/U_{\sf n}$. 
We then have an exact sequence $1 \rightarrow K^*/E_{\sf n} \rightarrow I_{\sf
n} \rightarrow {C}_{\sf n} \rightarrow 1$ (see loc. cit.) with $C_{\sf
n}$ finite and defined by means of this exact sequence. We can
consider the projective system of the $C_{{\sf n}p^r}$'s as $r$
varies and define $C_{{\sf n}p^{\infty}}$ to be the projective limit.
The character $\chi_{\wp}$ above maybe regarded naturally as a 
character of $C_{{\sf n}p^{\infty}}$
where $\sf n$ is the {\it conductor} of $\chi$.
 
Let $\chi$ be a Hecke
character. Since the kernel of the associated  homomorphism
$\chi_0^{}$ is open, there is a fractional ideal $a$, a
number field $L$  and a character $I_a\longrightarrow L^*$
induced by~$\chi$. Viewing $K^*$ as an algebraic torus, the
pull back to $K^*$ is algebraic and its kernel contains~$E$.
Conversely, a character $I_a\longrightarrow L^*$ whose
pull-back to
$K^*$ is algebraic automatically has $E$ in its kernel and
gives rise to a Hecke character. Just reverse the procedure
of going from $\chi$ to $\chi_0^{}$ above.

It will be convenient below to switch between ideal-theoretic and  
idele-theoretic viewpoints.
The strict ray class group ${\rm Cl}_{{\sf n}}$
sits inside the exact sequence 
$ 1 \rightarrow P_{{\sf n}} \rightarrow 
{\rm Id}_{{\sf n}} \rightarrow {\rm Cl}_{{\sf n}} \rightarrow 1$
where ${\rm Id}_{{\sf n}}$ is the group of fractional ideals
prime to ${\sf n}$ and $P_{\sf n}$ is the subgroup of
principal fractional ideals generated by $\gamma \in K$
with $\gamma$ totally positive and congruent to 1 mod ${\sf n}$.
Note that we have natural maps between exact sequences  
\begin{equation}
\begin{CD}
1 @> >> P_{{\sf n}} @> >>
{\rm Id}_{{\sf n}} @> >> {\rm Cl}_{{\sf n}}@>  >> 1
\\ 
&& @V VV @V VV @V VV 
\\
1@>  >> K^*/E_{\sf n}@>  >> I_{\sf
n}@>  >> {C}_{\sf n} @>  >> 1
\end{CD}
\end{equation}
 that
induces an isomorphism from ${\rm Cl}_{\sf n}$ to $C_{\sf n}$.

\subsection{Proof of Theorem \ref{main}}

Let the compatible system $\{\rho_{\wp}\}$ as in Theorem \ref{main}
have defect set $T$ and ramification set $S$.

By Lemma \ref{util} we may assume that $K$ is a Galois
extension of ${\bf Q}$ and is totally complex.
Let ${\sf m}_{\wp}$ be the ideal of ${\cal O}_K$ that is 
the prime to $p$ part of the Artin conductor of $\rho_{\wp}$
for places $\wp$ of $L$ not in $T$. We assume without loss of
generality that $L$ contains $K$.

{\it A priori} the fixed field of the kernel of $\rho_{\wp}$
is contained in the strict ray class field of $K$ of conductor ${\sf m}_{\wp}p^{\infty}$.
But as the image of $\rho_{\wp}$ has order prime to $p$, 
the fixed field of the kernel of $\rho_{\wp}$
is in fact contained in the ray class field of conductor ${\sf m}_{\wp}p$.
By class field theory we regard the data of the strictly compatible system
as the giving of homomorphisms $\rho_{\wp}:{\rm Cl}_{{\sf m}_{\wp}p}
\rightarrow {\bf F}_{\wp}^*$ where  ${\rm Cl}_{{\sf m}_{\wp}p}$
is the strict ray class group of $K$ of conductor ${\sf m}_{\wp}p$, and now
using the Artin map we see that the strict compatibility of the $\rho_{\wp}$'s
is expressed in terms of the images of prime ideals in 
${\rm Cl}_{{\sf m}_{\wp}p}$ under $\rho_{\wp}$. We will
for the most part only consider the restriction of  
$\rho_{\wp}:{\rm Cl}_{{\sf m}_{\wp}p} \rightarrow {\bf F}_{\wp}^*$
to the subgroup $P^{{\sf m}_{\wp}p}/P_{{\sf m}_{\wp}p}$ generated by the
principal ideals $P^{{\sf m}_{\wp}p}$ prime to ${\sf m}_{\wp}p$. 
We can also inflate this restriction and regard it as a homomorphism
of ${({\cal O}_K/{\sf m}_{\wp}p{\cal O}_K)}^*$ which factors through the
quotient by the image of the global units $E$.

Consider any 
integer $r$ of $K$ that generates a prime ideal that lies above a prime ideal
of ${\bf Q}$ which splits completely in $K$ and is unramified in $L$.
Abusing notation we will write
$\rho_{\wp}(r)$ for $\rho_{\wp}((r))$ especially when we are thinking of the homomorphism
of ${({\cal O}_K/{\sf m}_{\wp}p{\cal O}_K)}^*$ induced by $\rho_{\wp}$.
We see by the axioms of strictly compatible systems 
that there is an element 
$f_r \in L^*$ such that $\rho_{\wp}(r)$  is the reduction 
mod $\wp$  of $f_r$
that depends on $r$ but is independent of $\wp$,
and is $-f_r(0)$ using the notation of Definition \ref{compatible}. 
These considerations apply for all $r$ prime to $S$ and
for all primes $\wp$ not in $S \cup T$ and such that $p$, the residue
characteristic of $\wp$, does
not lie below $(r)$, and is prime to $f_r$. 

Choose a prime $\ell$ of ${\bf Q}$ that is coprime to the orders
of the multiplicative groups of the residue fields of primes in $S$
and is prime to the residue characteristics of the prime in $S$.
We claim that for almost all primes
$\wp$ of $L$ (in particular we exclude the places in $T$ and the
places which are ramified in $L$) 
the subgroup generated by the $\{\sigma(r)\}$'s in
$({\cal O}_L/\wp{\cal O}_L)^*$ contains the subgroup generated by $f_r^{t_{\wp}}$
(for some integer $t_{\wp}$ prime to $\ell$) 
in $({\cal O}_L/\wp{\cal O}_L)^*$, 
and where $\sigma$ runs through the distinct embeddings of $K$ in
$L$. The claim follows from considering the homomorphism of $({\cal
O}_K/{\sf m}_{\wp}p{\cal O}_K)^*$ ($\simeq({\cal
O}_K/{\sf m}_{\wp}{\cal O}_K)^* \times ({\cal
O}_K/p{\cal O}_K)^*$) induced by $\rho_{\wp}$. We are using here 
implicitly the fact that the map ${\cal O}_K/\wp' \rightarrow
{\cal O}_L/\wp$ is injective for any prime $\wp'$ of $K$ below $\wp$.
To prove the claim note that as
$\rho_{\wp}(r)$ is the reduction of $f_r$ mod $\wp$ and ${\sf
m}_{\wp}$ is only divisible by primes in $S$, 
the order of $f_r$ in ${\bf F}_{\wp}^*$, up to products of powers 
of primes which divide the orders
of the multiplicative groups of the residue fields of primes in $S$
and the residue characteristics of the primes in $S$, divides the order of 
$r \in ({\cal O}_K/p{\cal O}_K)^*$. The latter, by transport of structure,
is easily seen to be the l.c.m. of the orders of the
$\{\sigma(r)\}_{\sigma \in {\rm Gal}(K/{\bf Q})}$'s
in ${\bf F}_{\wp}^*$. 

We conclude from the claim and 
Proposition \ref{schoof} that $f_r^{t(r)}=\Pi_{\sigma}\sigma(r)^{m'_{r,\sigma}}$
for integers $t(r),m'_{r,\sigma}$ that a priori depend on $r$,
where $\sigma$ runs through the distinct embeddings of $K$ in $L$. 
Now as $r$ lies above a prime that is unramified in $L$, the $t(r)$'s divide all
the $m'_{r,\sigma}$'s. 
From this we conclude that we can write 
$f_r=\zeta_r\Pi_{\sigma}\sigma(r)^{m_{r,\sigma}}$ for integers
$m_{r,\sigma}$, where $\sigma$ runs through the distinct embeddings of $K$ in
$L$, and for some root of unity $\zeta_r$ that because of
our assumption that $L$ contains $K$ that is a Galois extension of
${\bf Q}$, and as $L$ contains only
finitely many roots of unity, has order bounded independently of $r$,
and thus divides a fixed integer say $n$.

Now consider two integers $r$ and $r'$ as above prime to $S$, that
generate prime ideals that 
lie above distinct primes ideals of ${\bf Q}$ that split completely in $K$ and
are unramified in $L$.
Then we see by as before that by considering the homomorphisms
of $({\cal O}_K/{\sf m}_{\wp}p)^*\rightarrow {\bf F}_{\wp}^*$ 
induced by $\{\rho_{\wp}\}$, and denoted by abuse of notation by the same symbol,
that $\rho_{\wp}(rr')$ is the reduction mod $\wp$ of
$f_rf_{r'}$. Thus
again as before we see that $(f_rf_{r'})=\Pi_{\sigma}\sigma((rr'))^{m_{rr',\sigma}}$
for some integers $m_{rr',\sigma}$. Thus we have an equality of ideals 
$$\Pi_{\sigma}\sigma((r))^{m_{r,\sigma}}\Pi_{\sigma}\sigma((r'))^{m_{r',\sigma}}=
\Pi_{\sigma}\sigma((rr'))^{m_{rr',\sigma}}$$ 
where $\sigma$ runs through the distinct embeddings of $K$ in
$L$. As $(r),(r')$ are split prime ideals in
$K$ lying above distinct primes of ${\bf Q}$, we  conclude
that in fact the $m_{r,\sigma}=m_{r',\sigma}$. 
For conjugates $\sigma(r)$ of $r$ apply the argument above now to $\sigma(r)r'$ and
thus we conclude that the $m_{r,\sigma}$'s are independent of $r$ and depend only
on $\sigma$. We denote the common value by $m_{\sigma}$. 
Note that at this point we have proved $f_r$ is integral at all places of $L$ not lying above the prime
of ${\bf Q}$ below $r$.

Recalling that
the $\rho_{\wp}$'s induce homomorphisms 
${({\cal O}_K/{\sf m}_{\wp}p{\cal O}_K)}^* \rightarrow {\bf
F}_{\wp}^*$ which factor through the
quotient by the image of the global units $E$, we see that the algebraic
character $\chi:K^* \rightarrow K^*$,
$\chi(x)=\Pi\sigma(x)^{nm_{\sigma}}$ contains the units $E$ of ${\cal
O}_K$ in its kernel. Then using the facts recalled
in  Section 4.1, and especially equation (1), 
the fact that follows from the Cebotarev density theorem
that the images of the split principal prime ideals coprime to ${\sf
m}_{\wp}p$
generate $P^{{\sf m}_{\wp}p}/P_{{\sf m}_{\wp}p}$,
and the axioms of strictly compatible systems, we see that 
for some Hecke character $\chi'$ of
infinity type $(nm_{\sigma})$, $\rho_{\wp}^n \otimes \tilde{\chi}_{\wp}^{'-1}$
factors through the Galois group of a fixed finite extension of $K$:
in fact this fixed finite extension can be taken to be
the Hilbert class field of $K$. Consequently by Lemma \ref{util} the 
strictly compatible system 
$\{\rho_{\wp}^n\}$ arises from a Hecke character
of infinity type $(nm_{\sigma})$. 

Hence by inspection as $\{\chi_{\wp}^n\}$ has trivial prime to $p$ conductor, 
we deduce that the compatible system $\{\rho_{\wp}\}$ is such
that the exponents of the primes dividing the prime to $p$ part of
the conductor of $\rho_{\wp}$ are bounded independently of $\wp$ and
thus as the exceptional set $S$ is finite the strictly compatible system
$\{\rho_{\wp}\}$ has {\it bounded conductor} in the sense of Definition
\ref{compatible}.

Thus we can take the ${\sf m}_{\wp}$'s to be independent
of $\wp$, and denote the common ideal by $\sf m$. Now essentially we
have
to repeat the argument above. Consider
principal prime ideals $(r)$ with $r$ congruent to 1 mod $\sf m$
that lie above primes that split completely in $K$,
and repeat the argument above to get that this time 
$f_r$, which we know {\it a priori} is in $K$ using 
the Frobenius density theorem, 
is $=\Pi_{\sigma}\sigma(r)^{m_{\sigma}}$ for the same $m_{\sigma}$'s as
above, using the last sentence of Proposition \ref{schoof}.
Now observe that the algebraic character $K^* \rightarrow K^*$ that sends $x$
to $\Pi_{\sigma}\sigma(x)^{m_{\sigma}}$ is trivial on units congruent
to 1 mod $\sf m$. Then using the facts recalled in Section 4.1,
and that such $r$'s, prime to $p$, project surjectively
to $({\cal O}_K/p{\cal O}_K)^*$ for almost all $p$,
we deduce by the axioms of strictly
compatible systems, 
that there is a Hecke character $\chi$ of
infinity type $(m_{\sigma})$, such that $\rho_{\wp} \otimes \tilde{\chi}_{\wp}^{-1}$
factors through the Galois group of a fixed finite extension of $K$.
The proof of the theorem is now complete by appealing to Lemma \ref{util}.

\subsection{Some remarks}

1. The proof of Theorem \ref{main} follows the 
general lines of the method of [K] that dealt with the case $K={\bf Q}$.
The presentation of the proof in [K] is 
inaccurate when $S$ is non-empty, as the second line of the
proof is unjustified (we thank N.~Fakhruddin for pointing this out):
nevertheless the proof of loc. cit. 
can be modified without much difficulty to work in the 
general case considered there of compatible systems of bounded
conductor. As compared to [K], the substantive improvements
made in this paper as far as the results about one-dimensional systems
are concerned, are that we no longer assume $K={\bf Q}$, we 
no longer have a bounded conductor hypothesis, and the the proof is
simplified to the extent that we no longer appeal to results 
towards Artin's conjecture on primitive roots, albeit we will need to use
results towards Artin's conjecture in the next section when studying
abelian semisimple compatible systems.

2. In [S2] a similar theorem was proved assuming that ``inertial weights'' (see Section 1.7 of [S3]) 
of the $\rho_{\wp}$'s were bounded independently of $\wp$.
In our proof the fact that inertial weights are bounded is proved to be a consequence
of the defining properties of a compatible system.

3. The proof works even if we allow the set $T$ in Definition
\ref{compatible} to be a set of places of density 0.

4.  The proof also gives that one-dimensional weakly compatible
mod $\wp$ systems are (weakly) equivalent in a natural sense to 
strictly compatible systems that arise from a 
uniquely determined Hecke character.   

5. Although our conjectures seem inaccessible at the moment in the
higher dimensional situation, it will be of interest
to prove some more accessible ``independence of $p$'' results for the images
of $\rho_{\wp}(G_K)$ analogous to the case of
compatible $\wp$-adic systems studied in [LP]. For instance one
might expect that for 2-dimensional strictly compatible systems either
for a set of primes $\wp$ of density 1, ${\rm im}(\rho_{\wp})$
has an abelian subgroup of bounded index,  or
for a set of primes $\wp$ of density 1, ${\rm im}(\rho_{\wp})$
contains $SL_2({\bf F}_p)$. This will be a necessary step
in studying the adelic images of compatible systems of $p$-adic
representations in the abstract case studied in [LP]. 

\section{Reciprocity for abelian semisimple compatible systems}

As a first step in addressing the conjectures of the introduction
for dimensions greater than 1, 
it is of interest to generalise Theorem \ref{main} to describe
abelian semisimple compatible mod $\wp$ systems of $G_K$ of arbitrary
dimensions. Such systems which are integral can be easily classified by using results in the 1-dimensional
case. We do not have satisfactory answers for rational abelian compatible systems. The problem in
generalising the proof of Theorem \ref{main} given above is
that we cannot use the argument given there to conclude that the
$m_{r,\sigma}$'s are independent of $r$. By using results towards Artin's conjecture
the case when $K={\bf Q}$ or more generally when 
$K$ is an abelian extension of ${\bf Q}$ can be treated. Below we state what can be proved and
only briefly indicate the arguments highlighting the novel features which arise
in the case of higher dimensional abelian systems. 

\subsection{Integral compatible systems}

We begin by indicating how a {\it d\'evissage} argument 
reduces the understanding of $L$-integral
abelian compatible systems to understanding the 1-dimensional case.

\begin{theorem}
Let $\{\rho_{\wp}\}$ be a $n$-dimensional 
compatible abelian $L$-integral semisimple
of $G_K$ of bounded conductor.
Then it arises from a direct sum of Hecke characters.
\end{theorem}

\noindent{\bf Sketch of proof:} For split principal prime ideals $(r)$ of $K$ with $r$ congruent to 1 
mod $\sf m$ ($\sf m$ is divisble by the prime to $p$ part of the conductor of $\rho_{\wp}$) we get as before that 
the roots of $f_r(X)$ are of the form
$\Pi_{\sigma}\sigma(r)^{m_{i,r,\sigma}}$, $i=1,\cdots,n$,
with the exponents non-negative by the integrality hypothesis. 

Looking at the compatible system
of representations of $G_K$ given by $\{{\rm det}(\rho_{\wp})\}$ we get a compatible system of
one-dimensional representations of $G_K$
of bounded conductor that by 
Theorem \ref{main} arises from a Hecke character. From this it is easy to
conclude that $\sum_i\sum_{\sigma}m_{i,r,\sigma}$ 
is bounded independently of $r$
and using the non-negativity of the exponents $m_{i,r,\sigma}$ we get that
the ``infinity types'' $m_{i,r,\sigma}$'s are bounded independently of $r$. Thus there are only
finitely many possibilities for the $m_{i,r,\sigma}$'s as $i,r,\sigma$ vary
and let $N$ be the sum of all these finitely many possibilities.
Let $(\alpha)$ be a split prime ideal of $K$ 
and choose a prime $p$ large enough such that whenever
$\Pi_\sigma \sigma(\alpha)^{m_{\sigma}}-1$ is not coprime to $p$
and $|m_{\sigma}| \leq N$ then all the $m_{\sigma}$'s are 0.
Consider integral elements $\beta$ in $K$ such that $\beta$ is congruent to $\alpha$ mod $p$,
and congruent to 1 mod $\sf m$,
and $\beta$ generates a split prime ideal.
We claim that the roots of $f_{\alpha}(X)$ and
$f_{\beta}(X)$ have the same ``infinity type''. This is because
the roots of these polynomials, which are of the form
$\Pi_{\sigma \in {\rm Gal}(K/{\bf Q})}\sigma(\alpha)^{m_{i,\alpha,\sigma}}$
and $\Pi_{\sigma \in {\rm Gal}(K/{\bf Q})}\sigma(\beta)^{m_{i,\beta,\sigma}}$
with the latter congruent to
$\Pi_{\sigma \in {\rm Gal}(K/{\bf Q})}
\sigma(\alpha)^{m_{i,\beta,\sigma}}$ mod $p$ by choice of $\beta$, 
are congruent mod $\wp$ under some ordering, for $\wp$ a prime above $p$. 
From this and the fact that $p$ was chosen so that, whenever
$\Pi_\sigma \sigma(\alpha)^{m_{\sigma}}-1$ is not coprime to $p$
and $|m_{\sigma}| \leq N$, then all the $m_{\sigma}$'s are 0,
the claim follows.
Such elements $\beta$ surject onto $({\cal O}_K/p'{\cal O}_K)^*$ 
for almost all primes $p'$ of ${\bf Z}$. From this we conclude that
$\{\rho_{\wp}\}$ arises from a direct sum of Hecke characters of infinity types that
can be read off from $f_{\alpha}(X)$.

\vspace{3mm}

\noindent{\bf Remark:} The proof extends to the case when we do not assume bounded
conductor hypothesis on remarking that ramification
indices of primes in ${\bf Q}_M$, the union of all extensions of ${\bf Q}$
of a {\it fixed} degree say $M$, are finite and the number of roots
of unity in ${\bf Q}_M$ is finite. 

\subsection{$K={\bf Q}$}

\begin{theorem}
Let $\{\rho_{\wp}\}$ be a compatible abelian $L$-rational semisimple
system of $G_{\bf
Q}$ with finite defect and exceptional sets and of bounded conductor. Then 
it arises from a direct sum of Hecke characters.
\end{theorem}

\noindent{\bf Sketch of proof:} Using known results towards Artin's conjecture (cf., [M] and [HB]) we find a prime
$q$ that is 1 mod $\sf m$ (with $\sf m$ as in the preceding paragraph) a primitive root mod $p$ for
infinitely many primes $p$. By using Proposition \ref{schoof} as
before we conclude that the roots of the polynomial $f_q(X)$, that is
part of the defining data of the compatible system $\{\rho_{\wp}\}$,
are $\{q^{m_1},\cdots,q^{m_n}\}$ for some integers $m_i$. Then we see
by a pigeonhole argument, using the infinitely many primes $p$
for which $q$ is a primitive root, that there are Dirichlet characters $\varepsilon_i$ of
conductor dividing $\sf m$ so that there are infinitely many primes $\wp$
such that $\rho_{\wp}$ is the direct sum $\oplus_{i=1}^n \varepsilon_{i,\wp}\widetilde{\chi_p^{m_i}}$
with $\chi_p$ the $p$-adic cyclotomic character. This proves the theorem.

\vspace{3mm}

\noindent{\bf Remark:} Using that
the closure of the subgroup generated by $q$ in $\Pi_{\ell \in S}{\bf Z}_{\ell}^*$
has finite index, it is easy to remove the bounded conductor hypothesis.

\subsection{A version of Artin's conjecture and abelian extensions $K$}

The case of rational abelian compatible systems of $G_K$ with $K$ a general number field
is one we are unable to resolve satisfactorily. In this section we assume that $K$ is abelian over ${\bf Q}$ and
indicate an approach.

We begin by formulating an Artin-type conjecture on primitive roots that may be of independent interest:
we formulate a very weak version that suffices for the purposes
here.

\begin{conj}\label{weakest}
  Let $K$ be a finite Galois extension of $\bf Q$ and $\sf m$ any
  ideal of ${\cal O}_K$. Then there is a totally positive
  integer $a \in {\cal O}_K$, $\simeq 1$ mod $\sf m$, which generates
  a split prime ideal of ${\cal
  O}_K$, and infinitely many primes $p$ of $\bf Q$ such that the conjugates of $a$,
  $\{\sigma(a)\} _{\sigma \in {\rm Gal}(K/{\bf Q})}$,
  generate a subgroup of $({\cal O}_K/p{\cal O}_K)^*/{\cal O}_K^*$ whose index is
  bounded idependently of $p$.
\end{conj}

\begin{theorem}
 Assume Conjecture \ref{weakest} and $K$ is an abelian extension of
 ${\bf Q}$. A $K$-rational compatible system $\{\rho_{\wp}\}$ 
 of abelian semisimple mod $\wp$ representations of
 ${\rm Gal}(\overline{K}/K)$ with bounded conductor arises from a direct sum of Hecke characters.
\end{theorem}

\noindent{\bf Sketch of proof:} Assuming the conjecture,
we choose a split principal, prime ideal $(a)$ that is congruent to 1 mod $\sf m$ where
$\sf m$ is an ideal invariant under ${\rm Gal}(K/{\bf Q})$
and which is divisible by the prime to $\wp$ part of the
conductor of $\rho_{\wp}$ for almost all $\wp$, such that such that for infinitely many primes $p$ of $\bf Q$ 
the elements $\{\sigma(a)\} _{\sigma \in {\rm Gal}(K/{\bf Q})}$
  generate a subgroup of $({\cal O}_K/p{\cal O}_K)^*/{\cal O}_K^*$ whose index is
  bounded idependently of $p$. Arguing as before
we know that the roots of $f_a(X)$ are of
the form $\Pi_{\sigma \in {\rm Gal}(K/{\bf Q})}\sigma(a)^{m_{i,a,\sigma}}$. Further we know
that for each $i$, the algebraic character of $K^*$ of infinity type
$(m_{i,a,\sigma})$ kills a subgroup of finite index of the units
${\cal O}_K^*$. 

After having noticed this, the main point
of interest of the proof of this theorem is the following proposition which we state for
general Galois extensions $K$ of ${\bf Q}$.
The point of it
is to ``algebraise'' an automorphism $\sigma$ of $K$
by interpreting it as the Frobenius map mod infinitely many primes
using Cebotarev density theorem.

\begin{prop}\label{equi} Let $\tau$ be in the centre of ${\rm Gal}(K/{\bf Q})$
with $K$ a finite Galois extension of ${\bf Q}$.
For an element $a$ of ${\cal O}_K$ as above,
for the positive density of primes $p$ of $Q$ such that
${\rm Frob}_{\wp} \in {\rm Gal}(K/{\bf Q})$ is $\tau$ 
for (all) $\wp$ above $p$ we have that 
$\rho_{\wp}({\rm Frob}_{\tau(a)})=\tau(\rho_{\wp}({\rm Frob}_{a}))$.
Thus in particular the characteristic polynomial of 
$\rho_{\wp}(\tau(a))$ is $\tau(f_a(X))$, for
any $\tau$ in the centre of ${\rm Gal}(K/{\bf Q})$ and for all 
$\wp$ not in $S \cup T$.
\end{prop}

\noindent{\bf Proof:} 
For the primes $p$ in the statement, $\tau$ naturally acts
on the domain and range of $\rho_{\wp}$. The proposition follows 
from \begin{itemize} \item $\tau$ maps the image of $a$ in
${\rm Cl}_{{\sf m}p}$ to its $p$th power as $a$ is 
congruent to 1 mod $\sf m$, and $\tau$ induces the Frobenius 
map on the residue fields of all primes of $K$ above $p$, \item and
$\tau$ induces the $p$th power map 
on ${\bf F}_{\wp}$.\end{itemize}

After this the proof of the theorem follows well-rehearsed lines.
We use the property of $a$ that it satisfies the conjecture above and
the proposition to conclude that the compatible system arises from
a Hecke character whose infinity type can be read off from the roots
of $f_a(X)$.  Namely we first ``twist'' $\{\rho_{\wp}\}$ by the
compatible system that arises from the direct sum of the Hecke characters of infinity types that
are determined by the roots of $f_a(X)$.
Here to make sense of ``twisting'' we are using that $\{\rho_{\wp}\}$ is abelian and we ``order'' the direct sum of Hecke characters
when twisting to match with the ordering of the roots of $f_a(X)$. This twisted sytem has the property
that the characteristic polynomials attached to $\sigma(a)$ are powers
of $X-1$. This uses the fact that for abelian $K$ 
the algebraic characters of $K^*$ are ${\rm Gal}(K/{\bf Q})$-equivariant.
Now as $a$ was chosen to satisfy the conjecture above we are done
by the usual argument that now we have a abelian semsimple system
$\{\rho_{\wp}'\}$ such that
$\rho_{\wp}'$ 
factors through the Galois group of a fixed extension of $K$
for infinitely many $\wp$, namely those $\wp$ for which $\langle \sigma(a)
\rangle$ generates a subgroup of bounded index of $({\cal O}_K/p{\cal O}_K)^*/{\cal O}_K^*$.

\vspace{3mm}

\noindent{\bf Remark:} In this case we do not know how to remove the bounded
conductor hypothesis. 

\section{Acknowledgements} I thank
Najmuddin Fakhruddin for his interest, and for the many interesting conversations
we've had about compatible systems. I thank the referee 
for useful suggestions to improve the paper.

\section{References}



\noindent [C] Chevalley, C., {\it  Deux th\'eor\`emes
d'arithm\'etique}, J. Math. Soc. Japan (1951), 36--44.

\vspace{3mm}

\noindent [Ca] Carayol, H., {\it Formes modulaires et repr\'esentations 
galoisiennes avec valeurs
dans un anneau local complet}, in {\it $p$-adic monodromy 
and the Birch and Swinnerton-Dyer conjecture}, Contemp. Math. 165 (1994),
213--237.

\vspace{3mm}

\noindent [Clo] Clozel, L., {\it Sur la th\'eorie de Wiles et le
changement de base nonab\'elien}, IMRN (1995), no. 9, 437--444.

\vspace{3mm}

\noindent [CS] Corrales, C., and Schoof, R., {\it The support
problem and its elliptic analogue}, J. of Number Theory 64 (1997), 276--290.

\vspace{3mm}

\noindent [DS] Deligne, P., Serre, J-P., {\it Formes modulaires de 
poids 1}, Annales
de l'Ecole Normale Superieure 7 (1974), 507--530.

\vspace{3mm}

\noindent [FM] Fontaine, J-M., Mazur, B., {\it Geometric Galois 
representations}, Elliptic curves, modular forms, and Fermat's last 
theorem (Hong Kong, 1993), 41--78, Ser. Number Theory, I, 
Internat. Press, Cambridge, MA, 1995. 

\vspace{3mm}

\noindent [H] Henniart, G., {\it Repr\'esentations $\ell$-adiques 
ab\'eliennes}, in S\'eminaire  de Th\'eorie des
Nombres, Progress in Math. 22 (1982), 107--126, Birkhauser.

\vspace{3mm}

\noindent [HB] Heath-Brown, D. R., {\it Artin's conjecture 
for primitive roots}, Quart. J. Math. Oxford 37 (1986), no. 145, 27--38.

\vspace{3mm}

\noindent [K] Khare, C., {\it Compatible systems of mod $p$ Galois representations},
C. R. Acad. Sci., Paris (1996), t. 323, S\'erie I, 117--120.

\vspace{3mm}

\noindent [K1] Khare, C., {\it Remarks on mod $p$ forms of weight one}, International
Mathematical Research Notices, vol. 3 (1997), 127--133.
(Corrigendum: IMRN 1999, no. 18, pg. 1029.)

\vspace{3mm}

\noindent [KK] Khare, C., Kiming, I., {\it Mod $pq$ Galois
representations and Serre's conjecture}, to appear in Journal of Number
Theory, preprint available at \\
{\tt http://www.math.utah.edu/\~{ }shekhar/papers.html}

\vspace{3mm}

\noindent [LP] Larsen, M., Pink, R., {\it On $\ell$-independence 
of algebraic monodromy groups in compatible systems of
representations}, Invent. Math. 107 (1992), 603--636. 

\vspace{3mm}



\noindent [M] Murty, Ram, {\it Artin's conjecture for primitive
roots}, Math. Intelligencer, 10 (1988) 59-67.

\vspace{3mm}

\noindent [S] Serre, J-P., {\it Abelian $\ell$-adic representations
and elliptic curves}, Addison-Wesley, 1989.



\vspace{3mm}

\noindent [S1] Serre, J-P., {\it  Sur les
repr\'esentations modulaires de degr\'e 2 de ${\rm
Gal}(\overline{\bf Q}/{\bf Q})$}, Duke Math. J. \textbf{54} (1987),
179--230.

\vspace{3mm}

\noindent [S2] Serre, J-P., {\it R\'esum\'es des cours de 1970-71}, 
{\it Oeuvres}, Vol. II, no. 93.

\vspace{3mm}

\noindent [S3]  Serre, J-P., {\it  Propri\'{e}t\'{e}s galoisiennes des
points d'ordre finides courbes elliptiques}, Invent. Math. 15 (1972), 259--331.

\vspace{3mm}

\noindent [Sch] Schappacher, N., {\it Periods of Hecke characters},
SLNM 1301.

\vspace{3mm}

\noindent [Tani] Taniyama, Y., {\it $L$-functions of 
number fields and zeta functions of abelian varieties}, J. Math. Soc.
Japan 9 (1957), 330--366.

\vspace{3mm}



\noindent [W] Weil, A., {\it On a certain type of characters of the
idele-class group of an algebraic number field}, Collected Papers,
volume 2, 255--261.

\vspace{3mm}

\noindent School of Mathematics, TIFR, Homi Bhabha Road, Mumbai 400 005,
INDIA. e-mail: shekhar@math.tifr.res.in

\noindent Dept of Math, University of Utah,
155 S 1400 E, Salt Lake City, UT 84112, USA. e-mail: shekhar@math.utah.edu

\end{document}